\documentclass[twocolumn]{autart} %{ifacconf}

\usepackage[sort, numbers]{natbib}
\usepackage{graphicx}          % Include this line if your 
\usepackage{amsmath,amsfonts,amssymb}
\usepackage{mathtools, textcomp}
\usepackage{mathrsfs}
\usepackage{graphicx}
\usepackage[dvips]{epsfig}
\usepackage{epstopdf}
\usepackage{color}
\usepackage{dsfont}

\newcommand{\SO}{\ensuremath{\mathrm{SO(3)}}}

\newcommand{\Ta}{\ensuremath{\mathrm{T}}}
\newcommand{\T}{^{\mbox{\small T}}}
\newcommand{\so}{\ensuremath{\mathfrak{so}(3)}}
\newcommand{\SE}{\ensuremath{\mathrm{SE(3)}}}
\newcommand{\se}{\ensuremath{\mathfrak{se}(3)}}
\newcommand{\bR}{\ensuremath{\mathbb{R}}}
\newcommand{\bS}{\ensuremath{\mathbb{S}}}

\newcommand{\mrm}{\mathrm}

\newcommand{\diag}{\mbox{diag}}
\newcommand{\bbm}{\begin{bmatrix}}
\newcommand{\ebm}{\end{bmatrix}}
\newcommand{\matl}{\left[ \begin{array}}
\newcommand{\matr}{\end{array} \right]}
\newcommand{\be}{\begin{equation}}
\newcommand{\ee}{\end{equation}}
\newcommand{\bea}{\begin{eqnarray}}
\newcommand{\eea}{\end{eqnarray}}
\newcommand{\beas}{\begin{eqnarray*}}
\newcommand{\eeas}{\end{eqnarray*}}
\newcommand{\nn}{\nonumber}
\newcommand{\mC}{\mathcal{C}}
\newcommand{\cE}{\mathcal{E}}
\newcommand{\cL}{\mathcal{L}}
\newcommand{\cI}{\mathcal{I}}

\newcommand{\cM}{\mathcal{M}}

\newcommand{\cT}{\mathcal{T}}
\newcommand{\cU}{\mathcal{U}}

\newcommand{\di}{\mathrm{d}}

\newcommand{\tr}{\mathrm{trace}}

\newcommand{\lan}{\langle}
\newcommand{\ran}{\rangle}

\newcommand{\cS}{\mathcal{S}}

\newcommand{\ad}[1]{{\mathrm{ad}_{#1}}}          			%adjoint
\newcommand{\adast}[1]{{\mathrm{ad}_{#1}^\ast}}  			%adjoint
\newcommand{\Ad}[1]{{\mathrm{Ad}_{#1}}}  			%adjoint
\DeclareMathOperator{\expm}{{expm}}  			% exponential map

\newtheorem{theorem}{Theorem}[section]

\newtheorem{lemma}[theorem]{Lemma}

\newtheorem{remark}[theorem]{Remark}

\def\thesection{\arabic{section}}

\newcommand{\bi}{\begin{itemize}}
\newcommand{\ei}{\end{itemize}}

\DeclareMathAlphabet{\mathpzc}{OT1}{pzc}{m}{it}

\newcommand{\bJ}{\mathbb{J}}

\newcommand{\mpz}{\mathpzc}
\newcommand{\msg}{\mathsf{g}}
\newcommand{\msh}{\mathsf{h}}

\newcommand{\sS}{\mathsf{S}}
\newcommand{\sO}{\mathsf{O}}
\newcommand{\bD}{\mathbb{D}}

\begin{document}

\begin{frontmatter}

\title{Rigid Body Pose Estimation based on the Lagrange-d'Alembert Principle$^\star$}
%\thanksref{footnoteinfo}}
%\thanks[footnoteinfo]{This work was supported in part by the National Science Foundation.}

%\author[First]{Maziar Izadi}, \ 
%\author[First]{Amit K. Sanyal\thanksref{footnoteinfo}}, \ 
%\author[Second]{Randal W. Beard} \ and \ 
%\author[Third]{Timothy W. McLain}
%
%\thanks[footnoteinfo]{corresponding author. Tel: (575) 646-2580}
%\address[First]{Department of Mechanical and Aerospace Engineering, New Mexico State 
%University, Las Cruces, NM 88003 USA (e-mail: {\tt\small \{mi,asanyal\}@nmsu.edu}) }                                      
%\address[Second]{Electrical and Computer Engineering Department,
%Brigham Young University, Provo, UT 84602 USA (e-mail: {\tt\small beard@ee.byu.edu})}
%\address[Third]{Mechanical Engineering Department, Brigham
%Young University, Provo, UT 84602 USA (e-mail: {\tt\small tmclain@et.byu.edu})}

\thanks[footnoteinfo]{This paper was not presented at any IFAC
meeting. Corresponding author A.~K.~Sanyal. Tel. +1(575) 646-2580.}

\author[NMSU]{Maziar Izadi}\ead{mi@nmsu.edu}, % Add the
\author[SU]{Amit K. Sanyal}\ead{aksanyal@syr.edu}, % (ead) as shown
%\author[ECE_BYU]{Randal W. Beard}\ead{beard@ee.byu.edu}, % e-mail address
%\author[ME_BYU]{Timothy W. McLain}\ead{tmclain@et.byu.edu}

\address[NMSU]{Department of Mechanical and Aerospace Engineering, New Mexico State University, Las Cruces, NM 88003 USA.} % Please supply
\address[SU]{Department of Mechanical and Aerospace Engineering, Syracuse University, Syracuse, 
NY 13244 USA (previously with New Mexico State University).}
%\address[ME_BYU]{Mechanical Engineering Department, Brigham Young University, Provo, UT 84602 USA.}
     
%\begin{keyword}                           % Five to ten keywords,  
%Cicero; Catiline; orations.               % chosen from the IFAC 
%\end{keyword}                             % keyword list or with the 
                                          % help of the Automatica 
                                          % keyword wizard

\begin{abstract}                          % Abstract of not more than 250 words.
Stable estimation of rigid body pose and velocities from noisy measurements, without any 
knowledge of the dynamics model, is treated using the Lagrange-d'Alembert principle from 
variational mechanics.
With body-fixed optical and inertial sensor measurements, a Lagrangian is obtained as the 
difference between a kinetic energy-like term that is quadratic in velocity estimation error and 
the sum of two artificial potential functions; one obtained from a generalization of Wahba's 
function for attitude estimation and another which is quadratic in the position estimate error. An 
additional dissipation term that is linear in the velocity estimation error is introduced, and the 
Lagrange-d'Alembert principle is applied to the Lagrangian with this dissipation. A Lyapunov 
analysis shows that the state estimation scheme so obtained provides stable asymptotic 
convergence of state estimates to actual states in the absence of measurement noise, with an 
almost global domain of attraction. This estimation scheme is discretized for computer 
implementation using discrete variational mechanics, as a first order Lie group variational 
integrator. The continuous and discrete pose estimation schemes require optical measurements 
of at least three inertially fixed landmarks or beacons in order to estimate instantaneous pose. 
The discrete estimation scheme can also estimate velocities from such optical measurements.
Moreover, all states can be estimated during time periods when measurements of only two 
inertial vectors, the angular velocity vector, and one feature point position vector are available 
in body frame. In the presence of bounded measurement noise in the vector measurements, 
numerical simulations show that the estimated states converge to a bounded neighborhood 
of the actual states.
\end{abstract} % 216/254 words

\end{frontmatter}

\section{Introduction}\label{Sec1}
\vspace*{-.1in}\hspace{.1in}
Estimation of rigid body translational and rotational motion is indispensable for operations of 
spacecraft, unmanned aerial and underwater vehicles. Autonomous state estimation of a rigid 
body based on inertial vector measurement and visual feedback from stationary landmarks, in 
the absence of a dynamics model for the rigid body, is analyzed here. The estimation scheme 
proposed here can also be applied to {\em relative state} estimation with respect to moving 
objects \cite{Gaurav_ASR}. This estimation scheme can enhance the autonomy and reliability of 
unmanned vehicles in uncertain GPS-denied environments. Salient features of this estimation 
scheme are: (1) use of onboard optical and inertial sensors, with or without rate gyros, for 
autonomous navigation; (2) robustness to uncertainties and lack of knowledge of dynamics; 
(3) low computational complexity for easy implementation with onboard processors; (4) proven 
stability with large domain of attraction for state estimation errors;  
and (5) versatile enough to estimate motion with respect to stationary as well as moving objects. 
Robust state estimation of rigid bodies in the absence of complete knowledge of their dynamics, 
is required for their safe, reliable, and autonomous operations in poorly known conditions. In 
practice, the dynamics of a vehicle may not be perfectly known, especially when the vehicle is 
under the action of poorly known forces and moments. The scheme proposed here has a single, 
stable algorithm for the coupled translational and rotational motion of rigid bodies using 
onboard optical (which may include infra-red) and inertial sensors. This avoids the need for 
measurements from external sources, like GPS, which may not be available in indoor, underwater 
or cluttered environments \cite{leishman2014relative,miller2014tracking,amelin2014algorithm}.
 
%\vspace*{-.2in}\hspace{.1in}
Attitude estimators using unit quaternions for attitude representation may be {\em unstable 
in the sense of Lyapunov}, unless they identify antipodal quaternions with a single attitude. 
This is also the case for attitude control schemes based on continuous feedback of unit 
quaternions, as shown in~\cite{Bayadi2014almost,sanyal2009inertia, chaturvedi2011rigid}. One adverse 
consequence of these unstable estimation and control schemes is that they end up taking 
longer to converge compared with stable schemes under similar initial conditions and initial 
transient behavior. Continuous-time attitude observers and filtering schemes on 
$\SO$ and $\SE$ have been reported in, e.g., \cite{Khosravian2015Recursive,mabeda04,silvest08,sany,markSO3,mahapf08,bonmaro09,Vas1,Khosravian2015observers,rehbinder2003pose}. 
These estimators do not suffer from kinematic singularities like estimators using coordinate 
descriptions of attitude, and they do not suffer from unwinding as they do not use 
unit quaternions. The maximum-likelihood (minimum energy) filtering method of 
Mortensen~\cite{Mortensen} was recently applied to attitude estimation, resulting in a nonlinear 
attitude estimation scheme that seeks to minimize the stored ``energy" in measurement 
errors~\cite{aguhes06,Zamani2013minimum,ZamPhD}. 
This scheme is obtained by applying Hamilton-Jacobi-Bellman (HJB) theory~\cite{kirk} to 
the state space of attitude motion \cite{ZamPhD}. Since the HJB equation can only be 
approximately solved with increasingly unwieldy expressions for higher order 
approximations, the resulting filter is only ``near optimal" up to second order. Unlike filtering 
schemes that are based on approximate or ``near optimal" solutions of the HJB equation and 
do not have provable stability, the estimation scheme obtained here can be solved 
exactly, and is shown to be almost globally asymptotically stable. Moreover, unlike filters 
based on Kalman filtering, the estimator proposed here does not presume any knowledge 
of the statistics of the initial state estimate or the sensor noise. Indeed, for 
vector measurements using optical sensors with limited field-of-view, the probability 
distribution of measurement noise needs to have compact support, unlike standard Gaussian 
noise processes that are commonly used to describe such noisy measurements.

%\vspace*{-.2in}\hspace{.1in}
The variational attitude estimator recently appeared in \cite{Automatica,ACC2015,CDC2015VPE}, where it was 
shown to be almost globally asymptotically stable. Some of the advantages of this scheme over 
some commonly used competing schemes are reported in \cite{ICRA2015}. This paper is the   
variational estimation framework to coupled rotational (attitude) and translational motion, as 
exhibited by maneuvering vehicles like UAVs. In such applications, designing separate state 
estimators for the translational and rotational motions may not be effective and may lead to poor 
navigation. For navigation and tracking the motion of such vehicles, the approach proposed here 
for robust and stable estimation of the coupled translational and rotational motion will be more 
effective than de-coupled estimation of translational and rotational motion states. Moreover, like 
other vision-inertial navigation schemes~\cite{shen2013vision,shen2013rotor}, the estimation 
scheme proposed here does not rely on GPS. However, unlike many other vision-inertial 
estimation schemes, the estimation scheme proposed here can
be implemented without any direct velocity measurements.
Since rate gyros are usually corrupted by high noise content and bias \cite{Good2013ECC}, 
such a velocity measurement-free scheme can result in fault tolerance in the case of faults with rate 
gyros. Additionally, this estimation scheme can be extended to relative pose estimation between 
vehicles from optical measurements, without direct communications or measurements of relative velocities.

%\vspace*{-.2in}\hspace{.1in}
The contents of this article are organized as follows. In Section \ref{Sec2}, the problem of motion 
estimation of a rigid body using onboard optical and inertial sensors is introduced. The 
measurement model is introduced and rigid body states are related to these measurements. Section \ref{Sec3} introduces artificial energy terms representing the measurement residuals 
corresponding to the rigid body state estimates. The Lagrange-d'Alembert principle is applied to the Lagrangian 
constructed from these energy terms with a Rayleigh dissipation term linear in the velocity 
measurement residual, to give the continuous time state estimator. Particular versions of 
this estimation scheme are provided for the cases when direct velocity measurements are not 
available and when only angular velocity is directly measured. Section \ref{Sec4} proves the 
stability of the resulting variational estimator. It is shown that, in the absence of measurement 
noise, state estimates converge to actual states asymptotically and the domain of 
attraction is an open dense subset of the state space. In Section \ref{Sec5}, the variational pose
estimator is discretized as a Lie group variational integrator, by applying the discrete 
Lagrange-d'Alembert principle to discretizations of the Lagrangian and the dissipation term. This 
estimator is simulated numerically in Section \ref{Sec6}, for two cases: the case 
where at least three beacons are measured at each time instant; and the under-determined case, 
where occasionally less than three beacons are observed. For these simulations, true states 
of an aerial vehicle are generated using a given dynamics model. Optical/inertial measurements 
are generated, assuming bounded noise in sensor readings. Using these 
measurements, state estimates are shown to converge to a neighborhood of actual states, for 
both cases simulated. Finally, Section \ref{Sec7} lists the  contributions and possible future 
extensions of the work presented in this paper.

\vspace*{-1mm}
\section{Navigation using Optical and Inertial Sensors} \label{Sec2}%Lidar and/or 
\vspace*{-.1in}\hspace{.1in}
Consider a vehicle in spatial (rotational and translational) motion. 
\begin{figure}[htb!]
	\centering
	\includegraphics[width=0.475\textwidth]{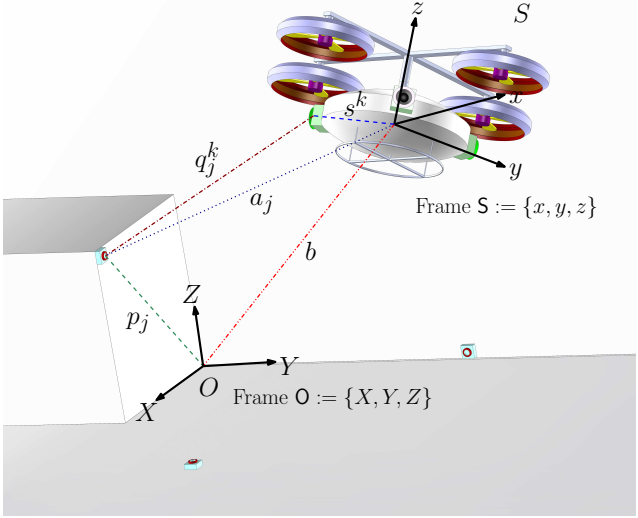}
    	\caption{Inertial landmarks on $O$ as observed from vehicle $S$ with optical  
	measurements.}
    	\label{Frames}
\end{figure}
Onboard estimation of the pose of the vehicle involves assigning a coordinate frame fixed to 
the vehicle body, and another coordinate frame fixed in the environment which takes the role 
of the inertial frame. Let $O$ denote the observed environment and $S$ denote the vehicle. 
Let $\mathsf{S}$ denote a coordinate frame fixed to $S$ and $\sO$ be a coordinate frame 
fixed to $O$, as shown in Fig. \ref{Frames}. Let $R\in\SO$ denote the rotation matrix from frame 
$\sS$ to frame $\sO$ and $b$ denote the position of origin of $\sS$ expressed in frame $\sO$.
%During initialization, the frame axes of 
%$\sO$ and $\sS$ could be chosen to be parallel.
The pose (transformation) from body fixed frame $\sS$ to inertial frame $\sO$ is then given by
\begin{align} 
&\msg=\bbm R \;\;\;& b\\ 0 \;\;\;& 1\ebm\in\SE.
\label{gDef} 
\end{align}
Consider vectors known in inertial frame $\sO$ measured by inertial sensors 
in the vehicle-fixed frame $\sS$; let $\beta$ be the number of such vectors. In addition, consider 
position vectors of a few stationary points in the inertial frame $\sO$ measured by optical (vision 
or lidar) sensors in the vehicle-fixed frame $\sS$. 
Velocities of the vehicle may be directly measured or can be estimated by linear 
filtering of the optical position vector measurements \cite{ACC2015}. 
Assume that these optical measurements are available for 
$\mpz j$ points at time $t$, whose positions are known in frame $\sO$ as $p_j$, $j\in\cI (t)$, 
where $\cI (t)$ denotes the index set of beacons observed at time $t$. Note that the observed 
stationary beacons or landmarks may vary over time due to the vehicle's motion. These points 
generate ${\mpz j\choose 2}$ unique relative position vectors, which are the vectors connecting 
any two of these landmarks. When two or more position vectors are optically measured, the number of vector measurements that
can be used to estimate attitude is ${\mpz j\choose 2}+\beta$. This number needs to be at least two 
(i.e., ${\mpz j\choose 2}+\beta\geq2$) at an instant, for the attitude to be uniquely determined at 
that instant. In other words, if at least two inertial vectors are measured at all instants (i.e., $\beta\geq2$), 
then beacon position measurements are not required for estimating attitude. However, at least one 
beacon or feature point position measurement is still required to estimate the position of the vehicle. 
Note that the use of two vector measurements for attitude determination was first proposed by the 
TRIAD algorithm in the 1960s \cite{TRIAD}.

\subsection{Pose Measurement Model} 
\vspace*{-.2in}\hspace{.1in}
Denote the position of an optical sensor and the unit vector from that sensor to an observed  
beacon in frame $\sS$ as $s^k\in\bR^3$ and $u^k\in\bS^2$, $k=1,\ldots,\mpz k$, respectively. 
Denote the relative position of the $j^{th}$ stationary beacon observed by the $k^{th}$ 
sensor expressed in frame $\sS$ as $q^k_j$. Thus, in the absence of measurement noise
\begin{align}
p_j=R(q^k_j+s^k)+b=Ra_j+b,\; j\in\mathcal I(t),
\label{FrameTrans}
\end{align}
where $a_j=q^k_j+s^k$, are positions of these points expressed in $\sS$. In practice, the $a_j$ 
are obtained from range measurements that have additive noise; we denote as $a_j^m$ the 
measured vectors. In the case of lidar range measurements, these are given by
\be a_j^m=(q^k_j)^m+s^k=(\varrho_j^k)^m u^k+s^k,\; j\in\mathcal I(t), \label{pim} \ee
where $(\varrho_j^k)^m$ is the measured range to the point by the $k^{th}$ sensor. 
The mean of the vectors $p_j$ and $a_j^m$ are denoted as $\bar{p}$ and $\bar{a}^m$ respectively, 
and satisfy
\begin{align}
\bar{a}^m= R\T (\bar p- b)+ \bar\varsigma, \label{barp}
%\bar{p}=R\bar{a}^m+b,\label{barp}
\end{align}
where $\bar{p}=\frac{1}{\mpz j}\sum\limits^\mpz j_{j=1}p_j$, $\bar{a}^m=\frac{1}{\mpz j}\sum
\limits^\mpz j_{j=1}a_j^m$ and $\bar\varsigma$ is the additive measurement noise 
obtained by averaging the measurement noise vectors for each of the $a_j$. Consider the 
${\mpz j\choose 2}$ relative position vectors from optical measurements, denoted as $d_j=
p_\lambda-p_\ell$ in frame $\sO$ and the corresponding vectors in frame $\sS$ as $l_j=
a_\lambda-a_\ell$, for $\lambda,\ell\in\cI (t)$, $\lambda\ne \ell$. The $\beta$ measured inertial 
vectors are included in the set of $d_j$, and their corresponding measured values expressed 
in frame $\sS$ are included in the set of $l_j$. If the total number of measured vectors (both 
optical and inertial), ${\mpz j\choose 2}+\beta=2$, then $l_3=l_1\times l_2$ is considered 
a third measured direction in frame $\sS$ with corresponding vector $d_3=d_1\times d_2$ in 
frame $\sO$. Therefore,
\begin{align}
d_j=Rl_j\Rightarrow D=RL,
\end{align}
where $D=[d_1\;\, \cdots\;\, d_n]$, $L=[l_1\;\, \cdots\;\, l_n]\in\bR^{3\times n}$ with $n=3$ if 
${\mpz j\choose 2}+\beta=2$ and $n={\mpz j\choose 2}+\beta$ if ${\mpz j\choose 2}+\beta>2$. 
Note that the matrix $D$ consists of vectors known in frame $\sO$. Denote the measured value 
of matrix $L$ in the presence of measurement noise as $L^m$. Then,
\begin{align}
L^m=R\T D+\mathscr{L},
\label{VecMeasMod}
\end{align}
where $\mathscr{L}\in\bR^{3\times n}$ consists of the additive noise in the vector measurements made
in the body frame $\sS$. 

\subsection{Velocities Measurement Model} 
Denote the angular and translational velocity of the rigid body expressed in body 
fixed frame $\sS$ by $\Omega$ and $\nu$, respectively. Therefore, one can write the kinematics 
of the rigid body as
\begin{align}
\dot{\Omega}=R\Omega^\times,\dot{b}=R\nu\Rightarrow\dot{\msg}= \msg\xi^\vee,
\label{Kinematics}
\end{align}
where $\xi= \bbm \Omega\\ \nu\ebm\in\bR^6$ and $\xi^\vee=\bbm \Omega^\times &\; \nu\\ 0 \;\;& 0\ebm$ and $(\cdot)^\times: \bR^3\to\so\subset\bR^{3\times 3}$ is the 
skew-symmetric cross-product operator that gives the vector space isomorphism between 
$\bR^3$ and $\so$:
\be
\mpz{x}^\times= \bbm \mpz{x}_1\\ \mpz{x}_2\\ \mpz{x}_3\ebm^\times= \bbm 0 & -\mpz{x}_3 & \mpz{x}_2\\ \mpz{x}_3 & 0 & -\mpz{x}_1\\ 
-\mpz{x}_2 & \mpz{x}_1 & 0\ebm.\label{timesmap}
\ee
For the general development of the motion estimation scheme, it is assumed that the velocities 
are directly measured. The estimator is then extended to cover the cases where: (i) only angular velocity is 
directly measured; and (ii) none of the velocities are directly measured.

\vspace*{-1mm}
\section{Dynamic Estimation of Motion from Proximity Measurements}\label{Sec3}
\vspace*{-.1in}\hspace{.1in}
In order to obtain state estimation schemes from measurements as outlined in Section 
\ref{Sec2} in continuous time, the Lagrange-d'Alembert principle is applied to an action functional of 
a Lagrangian of the state estimate errors, with a dissipation term linear 
in the velocities estimate error. This section presents the estimation scheme obtained 
using this approach. Denote the estimated pose and its kinematics as
\begin{align}
\hat\msg=\bbm \hat R & \;\;\; \hat{b}\\ 0 & \;\;\; 1\ebm\in\SE, \;\;\dot{\hat\msg}=\hat\msg\hat\xi^\vee,
\end{align}
where $\hat\xi$ is rigid body velocities estimate, with $\hat\msg_0$ as the initial pose estimate and the pose estimation 
error as
\begin{align}
\msh=\msg\hat\msg^{-1}=\bbm Q &\;\;\;\; b-Q\hat{b}\\ 0 & \;\;1\ebm=\bbm Q &\;\;\; x\\ 0 &\;\;\;1\ebm\in\SE,
\end{align}
where $Q=R\hat{R}\T$ is the attitude estimation error and $x=b-Q\hat{b}$. Then one obtains, in the case of perfect measurements,
\begin{align}\begin{split}
&\dot\msh = \msh\varphi^\vee,\, \mbox{ where }\, \varphi(\hat\msg,\xi^m,\hat\xi)=\bbm\omega\\\upsilon\ebm= \Ad{\hat\msg}\big(\xi^m- \hat\xi),
\end{split}\label{hdot}
\end{align}
where $\Ad{\mpz{g}}=\bbm \mpz{R}~~&~~0\\ \mpz{b}^\times\mpz{R}~~&~~\mpz{R}\ebm$ for $\mpz{g}=\bbm \mpz{R} &\;\;\; \mpz{b}\\ 0 & \;\;\;1\ebm$. The attitude and position estimation error dynamics are also in the form
\begin{align}
\dot{Q}=Q\omega^\times,\;\;\dot{x}=Q\upsilon.
\label{Qdot}
\end{align}

\vspace*{-1mm}
\subsection{Lagrangian from Measurement Residuals}
Consider the sum of rotational and translational measurement residuals between the 
measurements and estimated pose as a potential energy-like function. Defining the trace inner 
product on $\mathbb{R}^{n_1\times n_2}$ as
\begin{align}
\lan A_1,A_2\ran :=\tr(A_1 \T A_2),\label{tr_def}
\end{align}
the rotational potential function (Wahba's cost function~\cite{jo:wahba}) is expressed as
\begin{align}
\cU^0_r (\hat{\msg},L^m,D) &= \frac12\lan D -\hat R L^m , (D -\hat R L^m)W\ran,
\label{U0r}
\end{align}
where $W=\diag(w_j)\in\bR^{n\times n}$ is a positive diagonal matrix of weight factors for the 
measured $l_j^m$. Consider the translational potential function
\begin{align}
\cU_t (\hat \msg,\bar a^m,\bar p) &= \frac12 \kappa y\T y= \frac12\kappa \|\bar{p}-\hat{R}\bar{a}^m-\hat{b}\|^2,%\frac12\|\hat{R}(p_0^m-\hat{p}_0)\|^2.
\label{U0t} 
\end{align}
where $\bar{p}$ is defined by \eqref{barp}, $y\equiv y(\hat\msg,\bar{a}^m,\bar p)=\bar{p}-
\hat{R}\bar{a}^m-\hat{b}$ and $\kappa$ is a positive scalar. Therefore, the total potential function is 
defined as the sum of the generalization of \eqref{U0r} defined in~\cite{Automatica,ast_acc14}  
for attitude determination on $\SO$, and the translational energy \eqref{U0t} as
\begin{align}
\cU(\hat{\msg},L^m,D,\bar a^m,\bar p)&=\cU_r (\hat{\msg},L^m,D)+\cU_t(\hat{\msg},\bar a^m,\bar p)\nn\\
&= \Phi \big( \cU^0_r (\hat{\msg},L^m,D) \big)+\cU_t(\hat{\msg},\bar a^m,\bar p)\nn\\
&=\Phi \big(\frac12\lan D -\hat R L^m , (D -\hat R L^m)W\ran\big)\nn\\
&~~~~~~~+\frac12\kappa \|\bar{p}-\hat{R}\bar{a}^m-\hat{b}\|^2,\label{costU}
\end{align}
where $W$ is positive definite (not necessarily diagonal), and 
$\Phi: [0,\infty)\mapsto[0,\infty)$ is a $\mC^2$ function that satisfies $\Phi(0)=0$ and 
$\Phi'(\mpz x)>0$ for all $\mpz x\in[0,\infty)$. Furthermore, $\Phi'(\cdot)\leq\alpha(\cdot)$ where 
$\alpha(\cdot)$ is a Class-$\mathcal{K}$ function~\cite{khal} and $\Phi'(\cdot)$ denotes the 
derivative of $\Phi(\cdot)$ with respect to its argument. Because of these properties 
of the function $\Phi$, the critical points and their indices coincide for $\cU^0_r$ and 
$\cU_r$~\cite{Automatica}. Define the kinetic energy-like function: 
\be
%\cT (\hat\msg,\hat\xi,L^m,\mathcal{V}^m)=
\cT \Big(\varphi(\hat\msg,\xi^m,\hat\xi)\Big)= \frac12 \varphi(\hat\msg,\xi^m,\hat\xi)\T \bJ\varphi(\hat\msg,\xi^m,\hat\xi), 
\label{costT} \ee
where $\bJ\in\bR^{6\times 6}>0$ is an artificial inertia-like kernel matrix. Note that in contrast to 
rigid body inertia matrix, $\bJ$ is not subject to intrinsic physical constraints like the triangle 
inequality, which dictates that the sum of any two eigenvalues of the inertia matrix has to be 
larger than the third. Instead, $\bJ$ is a gain matrix that can be used to tune the estimator. For 
notational convenience, $\varphi(\hat\msg,\xi^m,\hat\xi)$ is denoted as $\varphi$ from 
now on; this quantity is the velocities estimation error in the absence of measurement 
noise. Now define the Lagrangian 
\be
%\cL (\hat \msg,\hat\xi,\mathcal{V}^m,L^m,D)= \cT (\hat\msg,\hat\xi,L^m,\mathcal{V}^m) %-\cU(\hat \msg,L^m,D),
\cL (\hat{\msg},L^m,D,\bar a^m,\bar p,\varphi)= \cT(\varphi) -\cU(\hat{\msg},L^m,D,\bar a^m,\bar p),
\label{contLag}\ee
and the corresponding action functional over an arbitrary time interval $[t_0,T]$ for $T>0$,
%\be \cS \big(\cL (\hat \msg,\hat\xi,\mathcal{V}^m,L^m,D)\big)= \int_{t_0}^T \cL 
%(\hat \msg,\hat\xi,\mathcal{V}^m,L^m,D) \di t, \,\label{action} \ee
\be \cS \big(\cL (\hat{\msg},L^m,D,\bar a^m,\bar p,\varphi)\big)= \int_{t_0}^T \cL 
(\hat{\msg},L^m,D,\bar a^m,\bar p,\varphi) \di t, \,\label{action} \ee
such that  $\dot{\hat \msg}= \hat \msg(\hat\xi)^\vee$. The following statement gives the form of 
the Lagrangian when perfect (noise-free) measurements are available, and derives the 
variational estimator for rigid body pose and velocities.
\begin{lemma}\label{NoiseFree}
In the absence of measurement noise, the Lagrangian is of the form
\begin{align}
\cL (\msh,D,\bar p,\varphi)=&\frac12 \varphi\T \bJ\varphi -\Phi \big(\lan I-Q, K\ran\big)-\frac12\kappa y\T y,\label{NFContLag}
\end{align}
where $K=DWD\T$ and $y\equiv y(\msh,\bar p)=Q\T x+(I-Q\T)\bar{p}$.
\end{lemma}
{\em Proof}: Suppose that all the measured states are noise free. Therefore, one can replace 
$L^m=L$, $\bar a^m=\bar a$ and $\xi^m=\xi$. The rotational potential function \eqref{U0r} can be 
replaced by
\begin{align}
\cU^0_r (\msh,D) &=\frac12\lan D-\hat{R}L^m, (D-\hat{R}L^m)W\ran\nn\\
&=\frac12\lan D-Q\T D, (D-Q\T D)W\ran\label{NFU0r}\\
&=\frac12\lan I-Q\T, (I-Q\T)DWD\T\ran=\lan I-Q, K\ran,\nn
\end{align}
since $\hat{R}E=Q\T D$ for the noise-free case. In addition,
\begin{align}
y(\msh, \bar p)&=\bar{p}-\hat{R}\bar{a}^m-\hat{b}=\bar{p}-\hat{R}\bar{a}-\hat{b}\\
&=\bar{p}-Q\T R\bar{a}-Q\T(b-x)=Q\T x+(I-Q\T)\bar{p}. \nn
\end{align}
The translational potential function in the absence of measurement noise can be expressed as
\begin{align}
\cU_t (\msh,\bar p) =\frac12\kappa y\T y.
\label{NFU0t}
\end{align}
Therefore, the total potential energy function is
\begin{align}
\cU(\msh,D,\bar p)&=\cU_r (\msh,D)+\cU_t(\msh,\bar p) \nn \\
&= \Phi \big( \cU^0_r (\msh,D) \big)+\cU_t(\msh,\bar p) \nn\\
&=\Phi \big(\lan I-Q, K\ran\big)+\frac12\kappa y\T y,
\label{NFcostU}
\end{align}
and the kinetic energy function is
\begin{align}
\cT(\varphi)= \frac12 \varphi\T \bJ\varphi.
\label{NFcostT}
\end{align}
Substituting \eqref{NFcostU} and \eqref{NFcostT} into:
\begin{align}
\cL (\msh,D,\bar{p},\varphi)&= \cT (\varphi) -\cU(\msh,D,\bar p)\nn\\
&=\cT (\varphi) -\Phi\big(\cU^0_r(\msh,D)\big)-\cU_t(\msh,\bar p),
\label{NFLag}
\end{align}
gives the Lagrangian \eqref{NFContLag} for the noise-free case.
\hfill\ensuremath{\square}

As in \cite{Automatica}, the positive definite weight matrix $W$ can be selected according to the 
following lemma:
\begin{lemma}\label{LemmaW}
Let $\mbox{rank}(D)=3$. Let the singular value decomposition of $D$ be given by
\begin{align}
D :&= U_D \Sigma_D V_D\T\, \mbox{ where }\,U_D\in\mathrm{O}(3),\ V_D\in\mathrm{O}(n),\nn\\
&\Sigma_D\in\mathrm{Diag}^+(3,n),
\end{align}
and $\mathrm{Diag}^+(n_1,n_2)$ is the vector space of $n_1\times n_2$ matrices with 
positive entries along the main diagonal and all other components zero. Let $\sigma_1, 
\sigma_2, \sigma_3$ denote the main diagonal entries of $\Sigma_D$. Further, let the positive 
definite weight matrix $W$ be given by 
\be W= V_D W_0 V_D\T\, \mbox{ where }\, W_0\in\mathrm{Diag}^+(n,n) \label{Wdec} \ee
and the first three diagonal entries of $W_0$ are given by
\be w_1= \frac{\varsigma_1}{\sigma_1^2},\; w_2=\frac{\varsigma_2}{\sigma_2^2},\; w_3=\frac{\varsigma_3}
{\sigma_3^2}\, \mbox{ where }\, \varsigma_1,\varsigma_2,\varsigma_3>0. \label{w123} \ee
Then, $K=DWD\T$ is positive definite and 
\be K= U_D\Delta U_D\T\, \mbox{ where }\, \Delta=
\diag(\varsigma_1,\varsigma_2,\varsigma_3), \ee
is its eigendecomposition. Moreover, if $\varsigma_\imath\ne \varsigma_\jmath\mbox{ for } \imath\ne \jmath$ and $\imath,\jmath\in\{1,2,3\}$, then $\lan I- Q,K\ran$ is a Morse function whose critical points are
\begin{align}
Q\in C_Q=\big\{ I, Q_1, Q_2, Q_3\big\} \mbox{ where } Q_\imath= 2 U_D I_\imath I_\imath\T U_D\T - I,
\label{C_Qdef}
\end{align}
and $I_\imath$ is the $\imath^{th}$ column vector of the identity $I\in\SO$.
\end{lemma}
The proof is presented in \cite{Automatica}.

\vspace*{-1mm}

\subsection{Variational Estimator for Pose and Velocities}
\vspace*{-.1in}\hspace{.1in}
The nonlinear variational estimator obtained by %expression \eqref{LagdAlem} 
applying the Lagrange-d'Alembert principle to the Lagrangian \eqref{contLag} with a dissipation 
term linear in the velocities estimation error, is given by the following statement.

\begin{theorem}\label{filterTHM}
The nonlinear variational  estimator for pose and velocities is given by
%that solves the variational problem \eqref{LagdAlem} 
\begin{align}
\begin{cases}
\bJ\dot{\varphi}&=\adast{\varphi}\bJ\varphi-Z(\hat{\msg},L^m,D,\bar{a}^m,\bar{p})-\bD\varphi,\vspace{.05in}\\
\hat{\xi}&=\xi^m-\Ad{\hat\msg^{-1}}\varphi,
\vspace{.05in}\\
\dot{\hat{\msg}}&=\hat{\msg}(\hat{\xi})^\vee,\label{ContFil}
\end{cases}
\end{align}
where $\adast{\zeta}=(\ad{\zeta})\T$ with $\ad{\zeta}$ defined by \eqref{ad_def}, and 
$Z(\hat{\msg},L^m,D,\bar{a}^m,\bar{p})$ is defined by
\begin{align}
\begin{split}
Z(\hat{\msg},L^m,D,&\bar{a}^m,\bar{p})=\label{Z}\\
&\bbm \Phi'\Big(\cU^0_r (\hat{\msg},L^m,D)\Big)S_\Gamma(\hat{R})+\kappa\bar{p}^\times y\\ \kappa y\ebm,\end{split}
\end{align}
where $\cU^0_r (\hat{\msg},L^m,D)$ is defined as \eqref{U0r}, $y\equiv y(\hat{\msg},\bar{a}^m,\bar{p})=\bar{p}-\hat{R}\bar{a}^m-\hat{b}$ and
\begin{align} S_\Gamma(\hat{R})&=\mrm{vex}\big(\Gamma\hat{R}\T-\hat{R}\Gamma\T\big)\nn\\
&=\mrm{vex}\big(DW(L^m)\T\hat{R}\T-\hat{R}L^mWD\T\big), \label{SLdef} \end{align}
$\Gamma= DW(L^m)\T$ and $\mrm{vex}(\cdot): \so\to\bR^3$ is the inverse of the $(\cdot)^\times$ map.
\end{theorem}
{\em Proof}: A Rayleigh dissipation term linear in the velocities of the form $\bD\varphi$ 
where $\bD\in\bR^{6\times 6} >0$ is used in addition to the Lagrangian \eqref{NFContLag}, and 
the Lagrange-d'Alembert principle from variational mechanics is applied to obtain the estimator on 
$\Ta\SE$.  {\em Reduced variations} with respect to $\msh$ and $\varphi$~\cite{blbook,marat} are 
applied, given by
\begin{align}
\delta\msh&= \msh\eta^\vee,\; \delta\varphi= \dot\eta+ \ad{\varphi}\eta,\label{variations}\\
\mbox{where }&\eta^\vee=\bbm \Sigma^\times & \rho\\ 0\;\; & 0\ebm\mbox{ and }\ad{\mpz{\zeta}}=\bbm \mpz w^\times~~ & 0\\ \mpz v^\times\;\; & \mpz{w}^\times\ebm,\label{ad_def}
\end{align}
for $\eta=\bbm \Sigma\\ \rho\ebm\in\bR^6$ and
$\zeta=\bbm \mpz w\\ \mpz v\ebm\in\bR^6$, with $\eta(t_0)=\eta(T)=0$. This leads to the expression:
%is applied by setting the first variation of the action functional
\begin{align}
\delta_{\msh,\varphi} \cS \big(\cL (\msh,D,\bar{p},\varphi)\big)=\int_{t_0}^T \eta\T
\bD\varphi \di t.\label{LagdAlem}
\end{align}
Note that the variations of the attitude and position estimation errors are of the form
\begin{align}
\delta Q=Q\Sigma^\times, \;\delta x=Q\rho,
\label{deltaQ}
\end{align}
respectively. Applying reduced variations to the rotational potential energy term 
\eqref{NFU0r}, one obtains
\begin{align}
\delta_{Q}\cU^0_r(\msh,D)&=\lan -Q\Sigma^\times, K\ran=\frac12\lan \Sigma^\times,KQ-Q\T K\ran
\nn\\
&=S_K\T(Q)\Sigma,
\label{deltaU0r}
\end{align}
where
\be S_K(Q)=\mrm{vex}\big(K Q-Q\T K\big). \label{SBdef} \ee
Taking first variation of the translational potential energy term \eqref{NFU0t} with respect to $Q$ 
and $x$ yields:
\begin{align}
\delta_{\msh}\cU_t(\msh,\bar{p})&=\kappa(\delta x+\delta Q\bar{p})\T\big\{x+(Q-I)\bar{p}\big\}\nn\\
&=\kappa\big(\rho\T y+\Sigma\T\bar{p}^\times y\big).
\label{deltaU0t}
\end{align}

Therefore, the first variation of the total potential energy \eqref{NFcostU} with respect to estimation errors is
\begin{align}
\delta_{\msh}\cU(\msh,D,\bar{p})=Z\T(\msh,D,\bar{p})\eta,
\end{align}
where $Z(\msh,D,\bar{p})$ is defined by
\begin{align}
Z&(\msh,D,\bar{p})=\label{NFZ}\\
&\bbm \Phi'\Big(\lan I-Q, K\ran\Big)S_K(Q)+\kappa\bar{p}^\times\big\{Q\T x+(I-Q\T)\bar{p}\big\} \\
\kappa \{Q\T x+(I-Q\T)\bar{p}\}\ebm.\nn
\end{align}
Taking the first variation of the kinetic energy term \eqref{NFcostT} with respect to $\varphi$ 
results in:
\begin{align}
\delta_{\varphi}\cT(\varphi)=\varphi\T\bJ\delta\varphi=\varphi\T\bJ(\dot{\eta}+\ad{\varphi}\eta),
\end{align}
applying the reduced variation for $\delta\varphi$ as given in \eqref{variations}. Therefore, the first 
variation of the action functional \eqref{action} is obtained as
\begin{align}
\delta&_{\msh,\varphi}\cS\big(\cL(\msh,D,\bar{p},\varphi)\big)\nn\\
&=\int_{t_0}^T\big\{\varphi\T\bJ(\dot{\eta}+\ad{\varphi}\eta)-\eta\T Z(\msh,D,\bar{p})\big\}\di t\nn\\
&=\int_{t_0}^T\eta\T\Big(\adast{\varphi}\bJ\varphi-Z(\msh,D,\bar{p})-\bJ\dot{\varphi}\Big)\di t+\varphi\T\bJ\eta|_{t_0}^T  \nn\\
&=\int_{t_0}^T\eta\T\Big(\adast{\varphi}\bJ\varphi-Z(\msh,D,\bar{p})-\bJ\dot{\varphi}\Big)\di t,
\label{actionVariation}
\end{align}
applying fixed endpoint variations with $\eta(t_0)=\eta(T)=0$. Substituting \eqref{actionVariation} in expression \eqref{LagdAlem} one obtains
\be \bJ\dot{\varphi}=\adast{\varphi}\bJ\varphi-Z(\msh,D,\bar{p})-\bD\varphi,\label{NFdotvarphi}\ee
where $Z(\msh,D,\bar{p})$ is defined by \eqref{NFZ}. In order to implement this estimator using the 
aforementioned measurements, substitute $Q\T D=\hat{R}L^m$. This changes the rotational 
potential energy formed by the estimation errors in attitude \eqref{NFU0r} to \eqref{U0r}.
%\begin{align}
%\cU^0_r (\msh,D) &=\lan I-Q,K\ran= \frac12\lan I-Q\T, (I-Q\T)DWD\T\ran\nn\\
%&=\frac12\lan (I-Q\T)D, (I-Q\T)DW\ran\nn\\
%&=\frac12\lan D-\hat{R}L^m, (D-\hat{R}L^m)W\ran.
%\label{nU0r}
%\end{align}
Equation \eqref{SBdef} is also reformulated as
\begin{align}
 S_K(Q)&=\mrm{vex}(DWD\T Q-Q\T DWD\T)\label{SBSL}\\
&=\mrm{vex}(DW(L^m)\T\hat{R}\T-\hat{R}(L^m)WD\T)=S_\Gamma(\hat{R}).\nn
\end{align}
Finally, the second row in the matrix $Z(\msh,D,\bar{p})$ is replaced by
\begin{align}
\kappa \{Q\T x+(I-Q\T)\bar{p}\}&=\kappa \{Q\T b-\hat{b}+\bar{p}-Q\T\bar{p}\}\nn\\
&=\kappa \{\hat{R}R\T(b-\bar{p})-\hat{b}+\bar{p}\}\nn\\
&=\kappa \{-\hat{R}\bar{a}^m-\hat{b}+\bar{p}\}.
\end{align}
Taking these changes into account, one could obtain the first of equations \eqref{ContFil} with 
$Z(\hat{\msg},L^m,D,\bar{a}^m,\bar{p})$ and $S_\Gamma(\hat{R})$ defined by \eqref{Z} and \eqref{SLdef}, respectively. 
Thus, the complete nonlinear estimator equations are given by \eqref{ContFil}.
\hfill\ensuremath{\square}

This is a fundamentally new idea of applying a principle from variational mechanics 
to obtain a state estimator, recently applied to rigid body attitude estimation 
in~\cite{Automatica}. This approach differs from the ``minimum-energy" approach to 
nonlinear estimation due to Mortensen~\cite{Mortensen} in some important ways. The 
minimum-energy approach applies Hamilton-Jacobi-Bellman (HJB) theory~\cite{kirk}, 
which can only be ``approximately solved." This approach was recently applied to state 
estimation of rigid body attitude motion in \cite{ZamPhD}. This HJB formulation can only be
approximately solved in practice, using a Riccati-like equation, to obtain a near-optimal filter 
that has no guarantees on stability. In the proposed approach, the time 
evolution of $(\hat \msg,\hat\xi)$ has the form of the dynamics of a rigid body with Rayleigh 
dissipation. This results in an estimator for the motion states $(\msg,\xi)$ that 
dissipates the ``energy" content in the estimation errors $(\mathsf{h},\varphi)= (\msg \hat\msg^{-1}, 
\Ad{\hat\msg}(\xi- \hat\xi))$ to provide guaranteed asymptotic stability in the case of perfect 
measurements~\cite{Automatica}. The differences between these two approaches were 
detailed in~\cite{ICRA2015}, for rigid body attitude estimation.

The proposed estimator combines certain desirable features of stochastic estimation and 
observer design approaches to state estimation for unmanned vehicles, when simultaneous 
inertial vector measurements and optical measurements of fixed beacons or landmarks are available. 
This nonlinear estimator is robust to measurement noise and does not require a dynamics model for 
the vehicle; instead, it estimates the dynamics of the vehicle given the measurement model 
in Section \ref{Sec2}. The variational pose estimator can also be interpreted as a low-pass stable filter (cf. \cite{Tayebi2011}). 
Indeed, one can connect the low-pass filter interpretation to the simple example of the natural dynamics 
of a mass-spring-damper system. This is a consequence of the fact that the mass-spring-damper system 
is a mechanical system with passive dissipation, evolving on a configuration space that is the vector 
space of real numbers, $\bR$. In fact, the equation of motion of this system can be obtained by 
application of the Lagrange-d'Alembert principle on the configuration space $\bR$. 
If this analogy or interpretation is extended to a system evolving on a Lie group as a 
configuration space, then the generalization of the mass-spring-damper system is a ``forced 
Euler-Poincar\'{e} system'' \cite{blbook,marat} with passive dissipation, as is obtained here. Explicit expressions for the vector of velocities $\xi^m$ can be obtained for two common cases 
when these velocities are not directly measured. These two cases are dealt with in the next subsection.

\subsection{Variational Estimator Implemented without Direct Velocity Measurements}\label{ContButterworth}

\vspace*{-1mm}
%\subsection{Linear Butterworth Filter for Measurements}\label{ContButterworth}
%\vspace*{-.1in}\hspace{.1in}
The velocity measurements in \eqref{ContFil} can be replaced by filtered velocity estimates 
obtained by linear filtering of optical and inertial measurements using, e.g., a second-order 
Butterworth filter. This is both useful and necessary when velocities are not directly measured. The 
filtered values $\xi^f$ are then used in place of $\xi^m$ to enhance the nonlinear estimator given by 
Theorem \ref{filterTHM}. Denote the measured vector quantity at time 
$t$ by $\mpz{z}^m$. A linear second-order filter of the form:
\begin{align}
\ddot{\mpz{z}}^f+ 2\mu\omega_n\dot{\mpz{z}}^f=\omega_n^2\big(\mpz{z}^m-\mpz{z}^f\big),
\label{ContBW}
\end{align}
is used, where $\omega_n$ is the natural (cutoff) frequency, $\mu$ is the damping ratio, and 
$\mpz{z}^f$ is the filtered value of $\mpz{z}^m$. Thereafter, %and $\mpz{z}^f_0=\mpz{z}^m(t_0)$. 
$\mpz{z}^f$ is used in place of $\mpz{z}^m$ in equations \eqref{ContFil}. 
%If $\mpz{z}^f$ is a filtered position 
%vector, then $\dot{\mpz{z}}^f$ is considered as the corresponding filtered velocity vector.

\subsubsection{Angular velocity is measured using rate gyros}
For the case that rate gyro measurements of angular velocities are available besides the 
$\mpz{j}$ feature point (or beacon) position measurements, the linear velocities of the rigid body can 
be calculated using each single position measurement by rewriting \eqref{pjdot} as
\begin{align}
\nu^f=(a_j^f)^\times\Omega^f-v_j^f,
\end{align}
for the $j^{th}$ point. Averaging the values of $\nu$ derived from all feature points gives a more reliable result. Therefore, the rigid body's filtered velocities are expressed in this case as
\begin{align}
\xi^f=\bbm \Omega^f\\ \frac{1}{\mpz{j}}\sum\limits_{j=1}^{\mpz{j}}(a_j^f)^\times\Omega^f-v_j^f \ebm.
\end{align}

\subsubsection{Translational and angular velocity measurements are not available}
In the case that both angular and translational velocity measurements are not available or accurate, 
rigid body velocities can be calculated in terms of the inertial and optical measurements. In order to do 
so, one can differentiate \eqref{FrameTrans} as follows
\begin{align}
&\dot{p}_j=R\Omega^\times a_j+R\dot{a}_j+\dot{b}=R\big(\Omega^\times a_j+\dot{a}_j+\nu\big)=0\nn\\
\Rightarrow&\dot{a}_j-a_j^\times\Omega+\nu=0\nn\\
\Rightarrow&v_j=\dot{a}_j=[a_j^\times\; -I]\xi=G(a_j)\xi,
\label{pjdot}
\end{align}
where $G(a_j)=[a_j^\times\; -I]$ has full row rank. From vision-based or Doppler lidar sensors, one can also 
measure the velocities of the observed points in frame $\sS$, denoted $v_i^m$. Here, velocity 
measurements as would be obtained from vision-based sensors is considered. 
The measurement model for the velocity is of the form
\be v_j^m=G(a_j)\xi+\vartheta_j,\ee
where $\vartheta_j\in\bR^3$ is the additive error in velocity measurement $v_j^m$. Instantaneous 
angular and translational velocity determination from such measurements is treated 
in~\cite{ast_acc14}. %For Doppler lidar measurements, $v_j^m=(\dot\varrho^k_j)^m u^k$, where 
%$\dot\varrho^k_j$ denotes the range-rate of point $j$ along line-of-sight $u^k$. 
Note that $v_j=\dot{a}_j$, for $j\in\mathcal I(t)$. As this kinematics indicates, the relative 
velocities of at least three beacons are needed to determine the vehicle's translational and 
angular velocities uniquely at each instant. However, when only one or two landmarks/beacons are 
measured, the estimator can propagate velocity estimates based on a least squares velocity 
determined from the available measurements. The rigid body velocities in both cases are obtained 
using the pseudo-inverse of $\mathds{G}(A^f)$:
\begin{align}
\mathds{G}(A^f)\xi^f&=\mathds{V}(V^f)\Rightarrow\xi^f=\mathds{G}^\ddag(A^f)\mathds{V}(V^f),\label{ximMore2}\\
\mbox{where }\;\mathds{G}(A^f)&=\bbm G(a^f_1)\\\vdots\\G(a^f_\mpz j)\ebm\;\mbox{and }\;\mathds{V}(V^f)=\bbm v^f_1\\\vdots\\v^f_\mpz j \ebm,\label{GVDef}
\end{align}
for $1,...,\mpz j\in\mathcal I(t)$. 
%, are matrices consisting of range rate measurements and a function 
%of range measurement for observed landmarks. 
When at least three beacons are measured, $\mathds{G}(A^f)$ is a full column rank matrix, 
and $\mathds{G}^\ddag(A^f)= \Big( \mathds{G}\T(A^f)\mathds{G} (A^f)\Big)^{-1} 
\mathds{G}\T (A^f)$ gives its pseudo-inverse. For the case that only one or two beacons are 
observed, $\mathds{G}(A^f)$ is a full row rank matrix, whose pseudo-inverse is given by 
$\mathds{G}^\ddag (A^f)= \mathds{G}\T (A^f)\Big( \mathds{G}(A^f)\mathds{G}\T 
(A^f)\Big)^{-1}$.

\vspace*{-1mm}
\section{Stability and Robustness of Estimator}\label{Sec4}
\vspace*{-.1in}\hspace{.1in}
The stability of the estimator (filter) given by Theorem \ref{filterTHM} is analyzed here.
The following result shows that this scheme is stable, with almost global convergence of 
the estimated states to the real states in the absence of measurement noise.
\begin{theorem} \label{thmFilt0}
Let the observed position vectors from optical measurements be bounded. Then, 
the estimator presented in Theorem \ref{filterTHM} is asymptotically stable at the estimation error 
state $(\msh,\varphi)=(I,0)$ in the absence of measurement noise. Further, the domain of 
attraction of $(\msh,\varphi)=(I,0)$ is a dense open subset of $\SE\times\bR^6$. 
\end{theorem}
%\begin{proof}
{\em Proof}: In the absence of measurement noise, $\hat R E=Q\T D$. Therefore, $\Phi\big(\cU^0_r
(\hat\msg,L^m,D)\big)=\Phi\big(\cU^0_r(\msh,D)\big)$ is a Morse function on $\SO$. The stability of this 
estimator can be shown using the following candidate Morse-Lyapunov function, which can be 
interpreted as the total energy function (equal in value to the Hamiltonian) corresponding to the 
Lagrangian \eqref{contLag}:
\begin{align}
V(\msh,&D,\bar{p},\varphi)=\cT (\varphi)+\cU(\msh,D,\bar{p})\label{lyapf} \\
&=\frac{1}{2}\varphi\T\bJ\varphi+\Phi \big(\lan I-Q, K\ran\big)+\frac12\kappa y\T y.\nn 
\end{align}
Note that $V(\msh,D,\bar{p},\varphi)\geq 0$ and $V(\msh,D,\bar{p},\varphi)=0$ if and only if $(\msh,\varphi)=(I,0)$. Therefore, $V(\msh,D,\bar{p},\varphi)$ is positive definite on 
$\SE\times\bR^6$. Using \eqref{Qdot}, one can derive the time derivative of \eqref{NFcostU} as
\begin{align}
\frac{\di}{\di t}\cU(\msh,D,&\bar{p})=\Phi'(\cU^0_r(\msh,D)\big)\lan-Q\omega^\times, K\ran+\kappa (\dot x+\dot Q\bar{p})\T (Qy)\nn\\
&=\Phi'(\cU^0_r(\msh,D)\big)\lan\omega^\times, -Q\T K\ran\nn\\
&~~~~~~~~~~~~~~~~~~~~~~~~~~~~~~~+\kappa (Q\upsilon+Q\omega^\times\bar{p})\T(Qy)\nn\\
&=\frac12\Phi'(\cU^0_r(\msh,D)\big)\lan\omega^\times, KQ-Q\T K\ran\nn\\
&~~~~~~~~~~~~~~~~~~~~~~~~~~~~~~~~~~~~~~~~~~+\kappa (\upsilon+\omega^\times\bar{p})\T y\nn\\
&=\Phi'(\cU^0_r(\msh,D)\big)S\T_K(Q)\omega+\kappa y\T \upsilon+\kappa(\bar{p}^\times y)\T \omega\nn\\
&=Z\T(\msh,D,\bar{p})\varphi,
\label{cU0dot} 
\end{align}
where $S_K(Q)$ is defined as \eqref{SBdef} and $Z(\msh,D,\bar{p})$ as \eqref{NFZ}. Therefore, the time derivative of the candidate Morse-Lyapunov function is
\begin{align}
\dot{V}(\msh,D,\bar{p}&,\varphi)=\varphi\T\bJ\dot\varphi+\varphi\T Z(\msh,D,\bar{p})\nn\\
&=\varphi\T\Big(\adast{\varphi}\bJ\varphi-Z(\msh,D,\bar{p})-\bD\varphi+Z(\msh,D,\bar{p})\Big)\nn\\
&=-\varphi\T\bD\varphi. \label{dlyapf}
\end{align}
noting that $\varphi\T\adast{\varphi}\bJ\varphi=0$. Hence, the derivative of the Morse-Lyapunov function is negative semi-definite. Note that the error dynamics for the pose estimate error $\msh$ is given by \eqref{hdot}, while the 
error dynamics for the velocities estimate error $\varphi$ is given by \eqref{NFdotvarphi}. Note that $D(t)$, as a function of time, is piecewise continuous and uniformly bounded. The first property 
(piecewise continuity) is naturally satisfied by $D(t)$, which is piecewise constant as the number and 
inertial positions of beacons (or feature points) observed by body-fixed optical sensors is piecewise 
continuous in time. The second property (uniform boundedness) is satisfied by $D(t)$ if the position 
vectors observed are bounded in $\bR^3$, as assumed in the statement. Therefore, the error dynamics 
for $(\msh,\varphi)$ is non-autonomous. Considering \eqref{lyapf} and \eqref{dlyapf}, and 
applying Theorem 8.4 in \cite{khal}, one can conclude that $\varphi\T \bD\varphi
\rightarrow 0$ as $t\rightarrow \infty$, which consequently implies $\varphi\rightarrow0$. 
Thus, the positive limit set for this system is contained in%where $\dot V(Q,\omega)=0$ is
\begin{align}
\cE = \dot{V}^{-1}(0)=\big\{(\msh,\varphi)\in\SE\times\se:\varphi\equiv0\big\}. 
\label{dotV0}
\end{align}
Substituting $\varphi\equiv 0$ in the first equation of the estimator \eqref{ContFil}, we obtain the 
positive limit set where $\dot V\equiv 0$ (or $\varphi\equiv 0$) as the set
\begin{align}
%\begin{split}
\mathscr{I} &= \big\{(\msh,\varphi)\in\SE\times\bR^6: Z(\msh,D,\bar{p})\equiv 0, 
\varphi\equiv0\big\} \label{invset}\\
&= \big\{(\msh,\varphi)\in\SE\times\bR^6: Q \in C_Q,\  Q\T x=0,\ \varphi\equiv0\big\},\nn
%\end{split} 
\end{align}
where $C_Q$ is defined by \eqref{C_Qdef}. Therefore, in the absence of measurement errors, all the 
solutions of this estimator converge 
asymptotically to the set $\mathscr{I}$. Define $\cU_r (Q):= \Phi \big(\lan I-Q, K\ran\big)$, 
which is the attitude measurement residual in the case of perfect measurements. 
Thus, the attitude estimate error converges to the set of critical points of $\cU_r (Q)$ in this 
intersection, and the position estimate error $x$ converges to zero. 
%
%lemma \eqref{lemma:MorseF}} yields that $\Phi(\lan I-\hat Q, K\ran)$ is a perfect Morse-function. Hence, it has a set of four non-degenerate critical points on $\SO$ as follows
%\begin{align}
%M_c=\big\{I,\diag(1,-1,-1),\diag(-1,1,-1),\diag(-1,-1,1)\big\}
%\label{MorseCPs}
%\end{align}
The unique global minimum of $\cU_r (Q)$ is at $Q=I$ (Lemma 2.1 in \cite{Automatica}), so 
this estimation error is asymptotically stable.\par
Now consider the set
\be \mathscr{C}= \mathscr{I}\setminus (I,0),  \label{othereqb} \ee
%= \big\{(\hat{Q},\hat{\omega})\in
%\SO\times\bR^3: \hat Q \in C_Q\setminus I,\ \hat{\omega}\equiv 0\big\},
which consists of all stationary states that the estimation errors may converge to, besides 
the desired estimation error state $(I,0)$. Note that all states in the stable manifold of a 
stationary state in $\mathscr{C}$ converge to this stationary state. 
%Therefore, we need to find the stable manifolds of the stationary states in $\mathscr{C}$. 
From the properties of the critical points $Q_\iota\in C_Q\setminus (I)$ of $\cU^0_r (Q)$, 
$(\iota=1,2,3)$ given in Lemma 2.1 of \cite{Automatica}, we see that the stationary points in 
$\mathscr{I}\setminus (I,0)=
\big\{ (\bbm Q_\iota\;\;&0\\0\;\;\;&1\ebm ,0) : Q_\iota\in C_Q\setminus (I)\big\}$ have stable manifolds 
whose dimensions depend on the index of $Q_\iota$. Since the velocities estimate error $\varphi$ 
converges globally to the zero vector, the dimension of the stable manifold $\cM^S_\iota$
of the critical points, i.e. $(\bbm Q_\iota\;\;&0\\0\;\;\;&1\ebm,0)\in\SE\times\bR^6$ is 
\be \dim (\cM^S_\iota) = 9+(3-\,\mbox{index of } Q_\iota)= 12- \,\mbox{index of } Q_\iota. \label{dimStabM} \ee
Therefore, the stable manifolds of $(\msh,\varphi)=(\bbm Q_\iota\;\;&0\\0\;\;\;&1\ebm,0)$ are nine-dimensional, 
ten-dimensional, or eleven-dimensional, depending on the index of $Q_\iota\in C_Q\setminus (I)$ 
according to \eqref{dimStabM}. Moreover, the value of the Lyapunov function $V(\msh,D,
\varphi)$ is non-decreasing (increasing when $(\msh,\varphi)\notin\mathscr{I}$) 
for trajectories on these manifolds when going backwards in time. This 
implies that the metric distance between error states $(\msh,\varphi)$ along these 
trajectories on the stable manifolds $\cM^S_\iota$ grows with the time separation between these 
states, and this property does not depend on the choice of the metric on $\SE\times\bR^6$. 
Therefore, these stable manifolds are embedded (closed) submanifolds of $\SE\times\bR^6$ 
and so is their union. Clearly, all states starting in the complement of this union, converge 
to the stable equilibrium $(\bbm Q_\iota\;\;&0\\0\;\;\;&1\ebm,0)=(I,0)$; therefore the domain of 
attraction of this equilibrium is
\[ \mbox{DOA}\{(I,0)\} = \SE\times\bR^6\setminus\big\{\cup_{\iota=1}^3 \cM^S_\iota\big\}, \]
which is a dense open subset of $\SE\times\bR^6$. 
\hfill\ensuremath{\square}

Therefore, the domain of attraction for the variational estimation scheme at $(\msh,\varphi)=(I,0)$ is 
almost global over the state space $\Ta\SE\simeq\SE\times\bR^6$, which is the best possible with 
continuous control and navigation schemes for systems evolving on a non-contractible state 
space~\cite{chaturvedi2011rigid,bo:miln}. In the presence of measurement noise with 
bounded frequencies and amplitudes, one can show that the expected values of the 
state estimates converge to a bounded neighborhood of the true states. The size of this 
neighborhood, which can be considered as a measure of the robustness of this estimation 
scheme, depends on the values of the estimator gains $\bJ$, $W$ and $\bD$. These estimator 
gains can be selected based on balancing the transient and steady-state behavior
of the estimator. 

\begin{remark}
In the special case that the weight matrix $W$ in Wahba's function is chosen as a piecewise time 
constant matrix according to Lemma \ref{LemmaW}, $K=DWD\T$ is a constant matrix for all time. 
Therefore, the RHS of \eqref{NFdotvarphi} is not explicitly dependent on time. This makes 
$(\msh,\varphi)$ an autonomous system and therefore the use of Theorem 8.4 of \cite{khal} is not 
required to prove asymptotic stability. One can apply LaSalle's invariance principle (Theorem 4.4 
in \cite{khal}) to prove the convergence of state estimates to the equilibrium $(I,0)$ in this case.
\end{remark}

\vspace*{-1mm}
\section{Discretization for Computer Implementation}\label{Sec5}
\vspace*{-.1in}\hspace{.1in}
For onboard computer implementation, the variational estimation scheme outlined above has to 
be discretized. This discretization is carried out in the framework of discrete geometric 
mechanics, and the resulting discrete-time estimator is in the form of a Lie group variational 
integrator (LGVI), as in~\cite{sany}. Since the estimation scheme proposed here is obtained from 
a variational principle of mechanics, it can be discretized by applying the discrete 
Lagrange-d'Alembert principle~\cite{marswest}. Consider an interval of time $[t_0, T]\in\bR^+$ 
separated into $N$ equal-length subintervals $[t_i,t_{i+1}]$ for $i=0,1,\ldots,N$, with $t_N=T$ 
and $t_{i+1}-t_i=\Delta t$ is the time step size. Let $(\hat \msg_i,\hat\xi_i)\in\SE\times\bR^6$ 
denote the discrete state estimate at time $t_i$, such that $(\hat \msg_i,\hat\xi_i)\approx 
(\hat \msg(t_i),\hat\xi(t_i))$ where $(\hat \msg(t),\hat\xi(t))$ is the exact solution of the 
continuous-time estimator at time $t\in [t_0, T]$. Let the values of the discrete-time measurements  
$\xi^m$, $\bar a^m$ and $L^m$ at time $t_i$ be denoted as $\xi^m_i$, $\bar a^m_i$ and $L^m_i$, 
respectively. Further, denote the corresponding values for the latter two quantities in inertial frame 
at time $t_i$ by $\bar p_i$ and $D_i$, respectively.
The term representing the energy content of the pose estimation error, given by 
\eqref{costU}, is discretized as
\begin{align}
\cU(\hat{\msg}_i,L^m_i,D_i,&\bar{a}^m_i,\bar{p}_i)=\cU_r (\hat{\msg}_i,L_i^m,D_i)+\cU_t(\hat{\msg}_i,
\bar{a}^m_i,\bar{p}_i)\nn\\
&= \Phi \big( \cU^0_r (\hat{\msg}_i,L_i^m,D_i) \big)+\cU_t(\hat{\msg}_i,\bar{a}^m_i,\bar{p}_i)\nn\\
&=\Phi \big(\frac12\lan D_i -\hat R_i L_i^m , (D_i -\hat R_i L_i^m)W_i\ran\big)\nn\\
&~~~~~~~~~~+\frac12\kappa\|\bar{p}_i-\hat{R}_i\bar{a}_i^m-\hat{b}_i\|^2,\label{DiscCostU}
\end{align}
where $W_i$ is the matrix of weight factors corresponding to $D_i$ at time $t_i$. The term 
encapsulating the energy in the velocities estimate error \eqref{costT}, is discretized as
\begin{align}
\cT \Big(\varphi(\hat\msg_i,\xi_i^m,\hat\xi_i)\Big)= \frac12 \varphi(\hat\msg_i,\xi_i^m,\hat\xi_i)\T \bJ\varphi(\hat\msg_i,\xi_i^m,\hat\xi_i), 
\end{align}
where $\bJ=\diag(J,M)$ and $M,J$ are positive definite matrices.

\begin{lemma}\label{NoiseFree}
In the absence of measurement noise, the discrete-time Lagrangian is of the form
\begin{align}
\cL (\msh_i,D_i,\bar{p}_i&,\varphi_i)=\frac12\lan\mathcal{J}\omega_i^\times,\omega_i^\times\ran+\frac12\lan M\upsilon_i,\upsilon_i\ran\label{DisLag}\\
&-\Phi \big(\lan I-Q_i,K_i\ran\big)-\frac12 \kappa y_i\T y_i,\nn
\end{align}
where $y_i\equiv y(\msh_i,\bar p_i)=Q_i\T x_i+(I-Q_i\T)\bar{p}_i$ and $\mathcal{J}$ is defined in terms of the 
matrix $J$ by $\mathcal{J}=\frac12\tr[J]I-J$.
\end{lemma}
A Lie group variational integrator (LGVI) introduced in \cite{Sanyal2011almost} is applied to the 
discrete-time Lagrangian \eqref{DisLag} to obtain the discrete-time filter.

\begin{theorem} \label{discfilter}
A first-order discretization of the estimator proposed in Theorem \ref{filterTHM} is given by
\begin{align}
(J\omega_i)^\times&=\frac{1}{\Delta t}(F_i\mathcal{J}-\mathcal{J}F_i\T),\label{LGVI_F}\\
(M+\Delta t\bD_t)\upsilon_{i+1}&=F_i\T M\upsilon_i\label{LGVI_upsilon}\\
&~~~~~~~~+\Delta t \kappa (\hat{b}_{i+1}+\hat{R}_{i+1}\bar{a}^m_{i+1}-\bar{p}_{i+1}),\nn\\
(J+\Delta t\bD_r)\omega_{i+1}&=F_i\T J\omega_i+\Delta t M\upsilon_{i+1}\times\upsilon_{i+1}\nn\\
+\Delta t&\kappa \bar{p}_{i+1}^\times (\hat{b}_{i+1}+\hat{R}_{i+1}\bar{a}^m_{i+1})\label{LGVI_omega}\\
-\Delta t&\Phi' \big( \cU^0_r (\hat{\msg}_{i+1},L_{i+1}^m,D_{i+1}) \big)S_{\Gamma_{i+1}}(\hat{R}_{i+1}),\nn\\
\hat\xi_i&=\xi^m_i-\Ad{\hat\msg_i^{-1}}\varphi_i,\label{LGVI_xihat}\\
\hat\msg_{i+1}&=\hat\msg_i\exp(\Delta t\hat\xi_i^\vee),\label{LGVI_ghat}
\end{align}
where $F_i\in\SO$, $\big(\hat\msg(t_0),\hat\xi(t_0)\big)=(\hat\msg_0,\hat\xi_0)$, %$\cU^0_r (\hat\msg_i,L_i^m,D_i)$ and $\cU^0_r (\msh,D)$ and
$\varphi_i=[\omega_i\T\;\upsilon_i\T]\T$, and $S_{\Gamma_i}(\hat R_i)$ is the value of  
$S_\Gamma(\hat R)$ at time $t_i$, with $S_\Gamma(\hat R)$ as defined by \eqref{SLdef}. 
%Further, $\varphi_i=[\omega_i\T\;\upsilon_i\T]\T$.
%, and $\mathds{G}$ and $\mathds{V}$ are matrix functions defined by \eqref{GVDef}.
\end{theorem}
%\;\;\xi^m_i= \frac{1}{\mpz j}\sum_{j=1}^\mpz j G^\ddag(a_{j_i}^m)v_{j_i}^m,

{\em Proof}: Consider first variations with fixed endpoints for the pose estimation errors in discrete 
time given by:
\begin{align}
\delta Q_i &=Q_i\Sigma_i^\times,\;\ \Sigma_0=\Sigma_N=0,\label{discQVar}\\
\delta x_i &=Q_i\rho_i,\;\;\;\; \rho_0=\rho_N=0, \label{discxVar}
\end{align}
where $\Sigma_i,\rho_i\in\bR^3$ are ``discrete variation vectors". It can be shown that for any 
$\omega\in\bR^3$ we have
\begin{align}
(J\omega)^\times=\omega^\times\mathcal{J}+\mathcal{J}\omega^\times.
\label{JomegaCross}
\end{align}
Discretizing \eqref{Qdot} assuming that the angular velocity estimation error is constant in the time interval $[t_i,t_{i+1}]$ with a constant time step size $\Delta t$, one gets
\begin{align}
Q_{i+1}=Q_iF_i, \;\ i\in\{0,1,2,\ldots,N-1\},
\label{AltDiscQ}
\end{align}
where $F_i\in\SO$ is given by
\begin{align}
F_i=\exp(\Delta t\omega_i^\times)\approx I+\Delta t\omega_i^\times.
\label{AltF}
\end{align}
The variation of $F_i$ can be derived from \eqref{AltDiscQ} and $\delta Q_i=Q_i\Sigma_i^\times$. Thus
\begin{align}
\delta F_i=-\Sigma_i^\times F_i+F_i\Sigma_{i+1}^\times.
\label{deltaF}
\end{align}
Using \eqref{JomegaCross} and \eqref{AltF}, one can enforce the skew-symmetry of 
$(J\omega_i)^\times$ by
\begin{align}
(J\omega_i)^\times&=\omega_i^\times\mathcal{J}+\mathcal{J}\omega_i^\times\approx\frac{1}{\Delta t}\Big((F_i-I)\mathcal{J}-\mathcal{J}(F_i\T-I)\Big)\nn\\
&=\frac{1}{\Delta t}(F_i\mathcal{J}-\mathcal{J}F_i\T).
\label{DiscJomegaCross}
\end{align}

From \eqref{hdot}, the continuous rate of change of the attitude estimation error is $\dot{x}=
Q\upsilon$, which can be approximated to first order in discrete-time as
\begin{align}
\frac{x_{i+1}-x_i}{\Delta t}\approx Q_i\upsilon_i\Rightarrow x_{i+1}=\Delta t Q_i\upsilon_i+x_i.
\label{Discx}
\end{align}
The first variation in $\upsilon_i$ is then calculated using \eqref{Discx} as
\begin{align}
\delta\upsilon_i&=\delta\Big(\frac{1}{\Delta t}Q_i\T(x_{i+1}-x_i)\Big)\nn\\
&=-\Sigma_i^\times\upsilon_i+\frac{1}{\Delta t}Q_i\T(\delta x_{i+1}-\delta x_i)\nn\\
&=-\Sigma_i^\times\upsilon_i+\frac{1}{\Delta t}F_i\rho_{i+1}-\frac{1}{\Delta t}\rho_i.
\label{deltanu}
\end{align}
The discrete Lagrangian \eqref{DisLag} can be rewritten as
\begin{align}
\cL(\msh_i,D_i,\bar{p}_i,&F_i,\upsilon_i)=\frac{1}{2\Delta t}\lan \mathcal{J}(F_i-I),(F_i-I)\ran\nn\\
&+\frac{\Delta t}{2}\lan M\upsilon_i,\upsilon_i\ran-\Delta t\Phi \big( \cU^0_r (\msh_i,D_i) \big)
\label{AltLag}\\
&-\frac{\Delta t}{2}\kappa(Q_iy_i)\T (Q_iy_i).\nn
\end{align}
The action functional \eqref{action} is replaced by the action sum 
\be \cS_d \big(\cL (\msh_i,D_i,\bar{p}_i,F_i,\upsilon_i)\big)= \Delta t
\sum_{i=0}^{N-1} \cL (\msh_i,D_i,\bar{p}_i,F_i,\upsilon_i). \label{actsum} \ee
Applying the discrete Lagrange-d'Alembert principle with two Rayleigh dissipation terms for 
angular and translational motions gives 
\begin{align}
\delta&\cS_d\big(\cL(\msh_i,D_i,\bar{p}_i,F_i,\upsilon_i)\big)\label{dLagdAl}\\
&~~~~~~~~~~~~~~~~~~~~~~~~~+\Delta t\sum_{i=0}^{N-1} \Big\{\lan\Sigma_i,\tau_i\ran+\lan 
\rho_i, f_i\ran\Big\}=0\nn\\
\Rightarrow&\sum_{i=0}^{N-1}\Bigg\{\frac{1}{\Delta t}\lan\delta F_i,\mathcal{J}(F_i-I)\ran+\Delta t\lan \delta\upsilon_i,M\upsilon_i\ran  \nn\\
-&\frac{\Delta t}{2}\Phi' \big( \cU^0_r (\msh_i,D_i) \big)\big\lan\Sigma_i^\times,S_{K_i}^\times(Q_i)\big\ran-\Delta t\kappa\lan \rho_i,y_i\ran\nn\\
-\Delta& t\kappa\lan \Sigma_i^\times,y_i\bar{p}_i\T\ran+\frac{\Delta t}{2}\lan\Sigma_i^\times,\tau_i^\times\ran+\Delta t\lan \rho_i,f_i\ran
\Bigg\}=0.\nn
\end{align}
As symmetric matrices are orthogonal to skew-symmetric matrices in the trace inner product, 
using \eqref{AltF} we can rewrite the first term in \eqref{AltLag} as
\begin{align}
\lan\delta F_i,\mathcal{J}(F_i-I)\ran&=\lan\Sigma_i^\times,\mathcal{J} F_i\T\ran-\lan\Sigma_{i+1}^\times,F_i\T\mathcal{J}\ran
\label{deltaFJ}\\
&=\frac12\lan\Sigma_i^\times,\mathcal{J} F_i\T\ran-\frac12\lan\Sigma_i^\times,F_i\mathcal{J}\ran\nn\\
&~~~~-\frac12\lan\Sigma_{i+1}^\times,F_i\T\mathcal{J}\ran+\frac12\lan\Sigma_{i+1}^\times,\mathcal{J}F_i\ran\nn\\
=-\frac{\Delta t}{2}\lan\Sigma_i^\times,&(J\omega_i)^\times\ran+\frac{\Delta t}{2}\lan\Sigma_{i+1}^\times,F_i\T(J\omega_i)^\times F_i\ran.\nn
\end{align}
Hence equation \eqref{dLagdAl} can be re-expressed as
\begin{align}
\sum_{i=0}^{N-1}&\Bigg\{ -\frac12\lan\Sigma_i^\times,(J\omega_i)^\times\ran+\frac12\lan\Sigma_{i+1}^\times,F_i\T(J\omega_i)^\times F_i\ran\nn\\
&-\frac{\Delta t}{2}\lan \Sigma_i^\times,(\upsilon_i\times M\upsilon_i)^\times\ran+\lan F_i\rho_{i+1},M\upsilon_i\ran\nn\\
&-\lan \rho_i,M\upsilon_i\ran-\frac{\Delta t}{2}\Phi' \big( \cU^0_r (\msh_i,D_i) \big)\big\lan\Sigma_i^\times,S_{K_i}^\times(Q_i)\big\ran\nn\\
&-\kappa\Delta t\big\lan \rho_i, y_i\big\ran-\frac{\kappa\Delta t}{2}\big\lan\Sigma^\times_i, (\bar{p}_i^\times y_i)^\times \big\ran\nn\\
&+\frac{\Delta t}{2}\lan\Sigma_i^\times,\tau_i^\times\ran+\Delta t\lan \rho_i,f_i\ran
\Bigg\}=0.
\label{AltLdA2}
\end{align}
Separating this equation into two (rotational and translational) parts leads to
\begin{align}
(M+\Delta t\bD_t)\upsilon_{i+1}&=F_i\T M\upsilon_i-\Delta t \kappa y_{i+1},\\
(J+\Delta t\bD_r)\omega_{i+1}&=F_i\T J\omega_i+\Delta t M\upsilon_{i+1}\times\upsilon_{i+1}\nn\\
&-\Delta t\kappa \bar{p}_{i+1}^\times y_{i+1}\\
&-\Delta t\Phi' \big( \cU^0_r (\msh_{i+1},D_{i+1}) \big)S_{K_{i+1}}(Q_{i+1}),\nn
\end{align}
using the identity $\mpz F\T\mpz w^\times \mpz F=(\mpz F\T\mpz w)^\times$ and by replacing 
$\tau_i=-\bD_r \omega_i$ and $f_i=-\bD_t \upsilon_i$, where $\bD_r$ and $\bD_t$ are positive 
definite matrices such that 
\[ \bD= \bbm \bD_r & 0\\ 0\; &\; \bD_t\ebm. \] 
In the presence of measurement noise, $Q_i\T D_i$ and $y_i$ %\equiv y(\hat{\msg}_i,\bar{a}^m_i,\bar{p}_i)$ 
are replaced by $\hat{R}_iL_i^m$ and $\bar{p}_i-\hat{b}_i-\hat{R}_i\bar{a}_i^m$, respectively. These 
give the discrete-time state estimator in the form of the Lie group variational integrator 
\eqref{LGVI_F}-\eqref{LGVI_ghat}. 
\hfill\ensuremath{\square}

Model-based discrete-time rigid body state estimators 
using LGVI schemes for attitude estimation were reported in \cite{sany,jgcd12}, but 
dynamics model-free state estimators using LGVIs have appeared only recently  
in~\cite{Automatica,ACC2015}.

\begin{remark}
In the absence of any direct velocity measurements or only angular velocity measurements, the expressions provided in Section \ref{ContButterworth} to calculate rigid body velocities are still valid in discrete-time. 
One can use the discrete-time variables introduced in this section in place of their continuous-time counterparts. The second-order Butterworth filter \eqref{ContBW} is discretized using the 
\textit{Newmark-$\beta$ Method} as follows:
\begin{align}
\begin{cases} 
\mpz{z}^f_{i+1}=\mpz{z}^f_i+\Delta t \dot{\mpz{z}}^f_i+\frac{\Delta t^2}{4}(\ddot{\mpz{z}}^f_i+\ddot{\mpz{z}}^f_{i+1})\\
\\
\dot{\mpz{z}}^f_{i+1}=\dot{\mpz{z}}^f_i+\frac{\Delta t}{2}(\ddot{\mpz{z}}^f_i+\ddot{\mpz{z}}^f_{i+1})
\end{cases}.
\label{Newmark}
\end{align}
Choosing $\omega_n=2$ and $\mu=\frac12$, this method gives the filtered positions and velocities as  follows:
\begin{align}
&\begin{Bmatrix}
\mpz{z}^f_{i+1}\\\dot{\mpz{z}}^f_{i+1}
\end{Bmatrix}
=\frac{1}{4+4\mu\omega_n\Delta t+\omega_n^2\Delta t^2}\\&
\begin{bmatrix}
4+4\mu\omega_n\Delta t-\omega_n^2\Delta t^2 \;&\; 4\Delta t\;&\; \omega_n^2\Delta t^2\\
-4\omega_n^2\Delta t \;&\; 4-4\mu\omega_n\Delta t-\omega_n^2\Delta t^2 \;&\; 2\omega_n^2\Delta t
\end{bmatrix}\nn\\
&~~~~~~~~~~~~~~~~~~~~~~~~~~~~~~~~~~~~~~~~~~~~~~~~~~~~~~~~~~\begin{Bmatrix}
\mpz{z}^f_i\\
\dot{\mpz{z}}^f_i\\
\mpz{z}_i^m+\mpz{z}_{i+1}^m
\end{Bmatrix}.\nn
\end{align}
where $\mpz{z}^m_i$ and $\mpz{z}^f_i$ are the corresponding value of quantities $\mpz{z}^m$ and $\mpz{z}^f$ at time instant $t_i$, respectively. As with the continuous time version, $\xi^m_i$ can be replaced 
with $\xi^f_i$ in the estimator equations. 
\end{remark}

\vspace*{-1mm}
\section{Numerical Simulations}\label{Sec6}
\begin{figure}
\begin{center}
\includegraphics[height=2.6in]{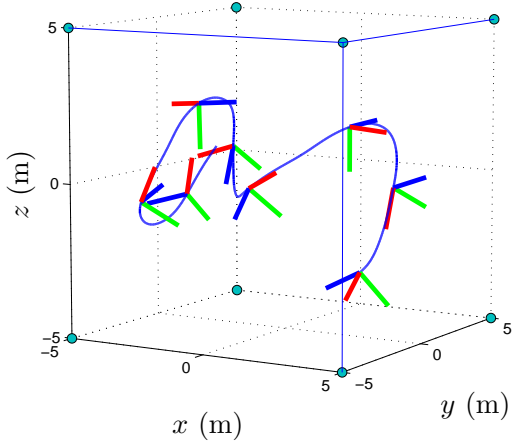}
\caption{Position and attitude trajectory of the simulated vehicle in 3D space.}  % width is 8.4 cm.
\label{Fig2}                                 % Size the figures 
\end{center}                                 % accordingly.
\end{figure}

\vspace*{-.1in}\hspace{.1in}
This section presents numerical simulation results for the discrete-time estimator obtained in 
Section \ref{Sec5}. In order to numerically simulate this estimator, simulated true states of an 
aerial vehicle flying in a room are produced using the kinematics and dynamics equations of a 
rigid body. The vehicle mass and moment of inertia are taken to be $m_v=420$ g and $J_v=
[51.2\;\;60.2\;\;59.6]\T$ g.m$^2$, respectively. The resultant external forces and torques applied 
on the vehicle are $\phi_v(t)=10^{-3}[10\cos(0.1t)\;\;2\sin(0.2t)\;\;-2\sin(0.5t)]\T$ N and 
$\tau_v(t)=10^{-6}\phi_v(t)$ N.m, respectively. The room is assumed to be a cubic space of size 
10m$\times$10m$\times$10m with the inertial frame origin at the center of this cube. The initial 
attitude and position of the vehicle are:
\begin{align}
R_0&=\expm_{\SO}\bigg(\Big(\frac{\pi}{4}\times[\frac{3}{7}\ -\frac{6}{7}\ \frac{2}{7}]\T\Big)^\times
\bigg),\nn\\
\mbox{and } b_0&=[2.5\ 0.5\ -3]\T\mbox{ m}.
\end{align}
This vehicle's initial angular and translational velocity respectively, are:
\begin{align}
\begin{split}
\Omega_0&=[0.2\;\; -0.05\;\; 0.1]\T\mbox{ rad/s},\\
\mbox{and } \nu_0&=[-0.05\;\;0.15\;\;0.03]\T\mbox{ m/s}.
\end{split}
\end{align}
The vehicle dynamics is simulated over a time interval of $T=150 \mbox{ s}$, with a time 
stepsize of $\Delta t=0.02 \mbox{ s}$. The trajectory of the vehicle over this time interval is 
depicted in Fig.~\ref{Fig2}.
The following two inertial directions, corresponding to nadir and Earth's magnetic field 
direction, are measured by the inertial sensors on the vehicle:
\begin{align}
d_1=[0\;\;0\;\;-1]\T,\;\;d_2=[0.1\;\;0.975\;\;-0.2]\T.
\end{align}
For optical measurements, eight beacons are located at the eight vertices of the cube, labeled 
1 to 8. The positions of these beacons are known in the inertial frame and their index (label) 
and relative positions are measured by optical sensors onboard the vehicle whenever the 
beacons come into the field of view of the sensors. Three identical cameras (optical sensors) 
and inertial sensors are assumed to be installed on the vehicle. The cameras are fixed to known 
positions on the vehicle, on a hypothetical horizontal plane passing through the vehicle, 120$^\circ$ 
apart from each other, as shown in Fig.~\ref{Frames}. All the camera readings contain random 
zero mean signals whose probability distributions are normalized bump functions with width of 
$0.001$m. The following are selected for the positive definite estimator gain matrices:
\begin{align}
%\begin{split}
J&=\diag\big([0.9\;\;0.6\;\;0.3]\big), \nn \\
M&=\diag\big([0.0608\;\;0.0486\;\;0.0365]\big), \\
%\end{split}
%\end{align}
%The rotational and translational motion's dissipation matrices are also selected as the follows:
%\begin{align}
%\begin{split}
\bD_r&=\diag\big([2.7\;2.2\;1.5]\big),\bD_t=\diag\big([0.1\;\;0.12\;\;0.14]\big). \nn
%\end{split}
\end{align}
$\Phi(\cdot)$ could be any $C^2$ function with the properties described in Section \ref{Sec3}, but is 
selected to be $\Phi(x)=x$ here. The initial state estimates have the following values:
\begin{align}
\begin{split}
\hat\msg_0&=I,\;\;\;\hat\Omega_0=[0.1\;\;0.45\;\;0.05]\T\mbox{ rad/s},\\
\mbox{ and }\hat\nu_0&=[2.05\;\;0.64\;\;1.29]\T\mbox{ m/s}.
\end{split}
\end{align}
The performance of the proposed estimator is presented for two different cases.

\vspace*{-1mm}
\subsection{CASE 1: At least three beacons are observed at each time instant}
\begin{figure}
\begin{center}
\includegraphics[height=2.4in]{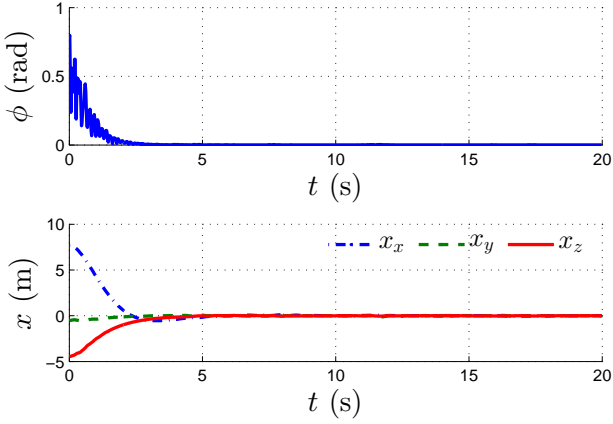}
\caption{Principal angle of the attitude and position estimation error for CASE 1.}  % width is 8.4 cm.
\label{Fig3}                                 % Size the figures 
\end{center}                                 % accordingly.
\end{figure}

\begin{figure}
\begin{center}
\includegraphics[height=2.4in]{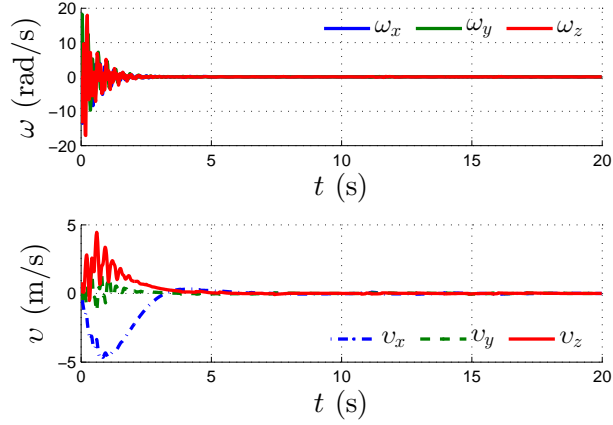}
\caption{Angular and translational velocity estimation error for CASE 1.}  % width is 8.4 cm.
\label{Fig4}                                 % Size the figures 
\end{center}                                 % accordingly.
\end{figure}

\vspace*{-.1in}\hspace{.1in}
Having three beacons measured at each time instant guarantees full determination of vehicle's 
translational and angular velocities instantaneously. A conic field of view (FOV) of 
2$\times$40$^\circ$ for cameras can satisfy this condition. The vehicle's velocity is calculated by 
\eqref{ximMore2} in this case. The discrete-time estimator \eqref{LGVI_F}-\eqref{LGVI_ghat} is 
simulated over a time interval of $T=20$ s with sampling interval $\Delta t=0.02$ s. At each time instant, 
\eqref{LGVI_F} is solved using the Newton-Raphson iterative method to find an approximation for 
$F_i$. Following this, the remaining equations (all explicit) are solved to generate the 
estimated states. The principal angle of the attitude estimation error and the position estimation  
error for CASE 1 are plotted in Fig.~\ref{Fig3}. Plots of the angular and translational velocity estimation  
errors are shown in Fig.~\ref{Fig4}.

\vspace*{-1mm}
\subsection{CASE 2: Less than three beacons are measured at some time instants}
\begin{figure}
\begin{center}
\includegraphics[height=2.4in]{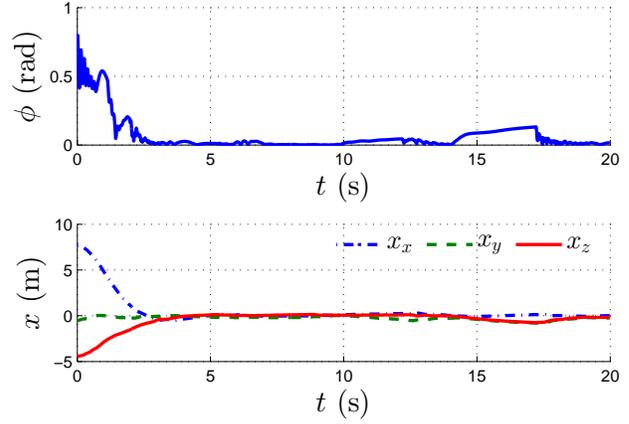}
\caption{Principal angle of the attitude and position estimation error for CASE 2.}  % width is 8.4 cm.
\label{Fig5}                                 % Size the figures 
\end{center}                                 % accordingly.
\end{figure}

\begin{figure}
\begin{center}
\includegraphics[height=2.4in]{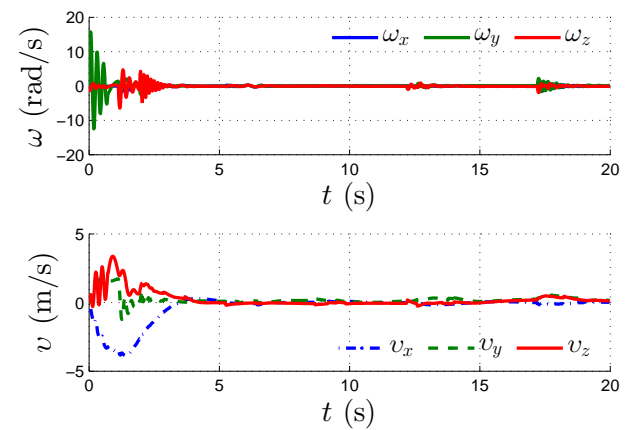}
\caption{Angular and translational velocity estimation error for CASE 2.}  % width is 8.4 cm.
\label{Fig6}                                 % Size the figures 
\end{center}                                 % accordingly.
\end{figure}

\vspace*{-.1in}\hspace{.1in}
To implement the variational estimator for the case that less than three optical 
measurements are available, the field of view of the cameras is decreased to limit the number of 
beacons observed. Assuming the cameras have conical fields of view of 2$\times$25$^\circ$, 
the minimum number of beacons observed instantaneously drops to 1 during 
the simulated time interval. The dynamics model for the aerial vehicle, simulated time duration, 
and sample rate are identical to CASE 1. Fig.~\ref{Fig5} depicts the principal angle of the 
attitude estimation error and the position estimation error for CASE 2, and Fig.~\ref{Fig6} shows 
the angular and translational velocity estimation errors. All estimation errors are shown to 
converge to a neighborhood of $(\msh,\varphi)=(I,0)$ in both cases, where the size of this 
neighborhood depends on the magnitude of measurement noise.

\vspace*{-1mm}
\section{Conclusion}\label{Sec7}
\vspace*{-.1in}\hspace{.1in}
This article proposes an estimator for rigid body pose and velocities, using optical and inertial 
measurements by sensors onboard the rigid body. The sensors are assumed to provide 
measurements in continuous-time or at a sufficiently high frequency, with bounded 
measurement noise. An artificial kinetic energy quadratic in rigid body velocity estimate errors is 
defined, as well as two fictitious potential energies: (1) a generalized Wahba's cost function for 
attitude estimation error in the form of a Morse function, and (2) a quadratic function of  
the vehicle's position estimate error. Applying the Lagrange-d'Alembert principle on a Lagrangian 
consisting of these energy-like terms and a dissipation term linear in velocities estimation error, an 
estimator is designed on the Lie group of rigid body motions. In the absence of measurement 
noise, this estimator is shown to be almost globally asymptotically stable, with estimates converging 
to actual states in a domain of attraction that is open and dense in the state space. The continuous 
estimator is discretized by applying the discrete Lagrange-d'Alembert principle on the discrete 
Lagrangian and dissipation terms linear in rotational and translational velocity estimation errors. 
In the presence of measurement noise, numerical simulations show that state estimates converge 
to a bounded neighborhood of the true states. Future extensions of this work include higher-order 
discretizations of the continuous-time filter given here and obtaining a stochastic 
interpretation of the variational pose estimator.

%\begin{ack}                               % Place acknowledgements
%\vspace*{-.1in}\hspace{0.085in}Support from the National Science Foundation through grant CMMI 1131643 is gratefully 
%acknowledged.
%\end{ack}

\end{document}